\documentclass{amsart}
\usepackage{amsmath,amsthm,amssymb,amscd,amsfonts}
\usepackage{enumerate}
\makeatletter
\@namedef{subjclassname@2010}{%
  \textup{2010} Mathematics Subject Classification}
\makeatother
\theoremstyle{definition}

\newcommand{\cA}{\mathcal A}

\newtheorem{theo}{Theorem}[section]

\theoremstyle{definition}

\theoremstyle{definition}

\numberwithin{equation}{section}

\title{ $\varphi$-contractibility and $\varphi$-Connes amenability coincide with some older notions}
\author[    A. Mahmoodi]{  Amin Mahmoodi}
\thanks{}

\address{Department of Mathematics, Central Tehran Branch,
 Islamic Azad University, Tehran, Iran}
\address{}

  \email{a\_mahmoodi@iauctb.ac.ir }
\begin{document}
\maketitle

\setcounter{section}{0}
\begin{abstract}
It is shown that various definitions of $\varphi$-Connes amenability
as introduced independently in \cite{Gh-Ja, mah, Sh-Am}, are just
rediscovering existing notions and presenting them in different
ways. It is also proved that even $\varphi$-contractibility as
defined in \cite{Sangani}, is equivalent to an older and simpler
concept.
\end{abstract}

\section{Introduction}
The precious and fertile notion of amenability was introduced by B.
E. Johnson in \cite{Joh}. A generalization of amenability depending
on homomorphisms was introduced and studied by E. Kaniuth, A. T. Lau
and J. Pym \cite{Lau}, and independently by M. S. Monfared
\cite{Sangani1}. For a Banach algebra $\mathfrak{A}$, we write
$\Delta(\mathfrak{A})$ for the set of all homomorphisms from
$\mathfrak{A}$ onto $ \mathbb{C}$. Let $\varphi \in
\Delta(\mathfrak{A})$. An element $m \in \mathfrak{A}^{**}$ is
called a \textit{right} [\textit{left}] $\varphi$-\textit{mean} if $
m(\varphi)= 1$ and $ m( f \cdot a) = \varphi(a) m(f)$ [$ m( a \cdot
f) = \varphi(a) m(f)$] for $a \in \mathfrak{A}$ and $f \in
\mathfrak{A}^*$. A Banach algebra is \textit{right} [\textit{left}]
$\varphi$-\textit{amenable} if it has a right [left] $\varphi$-mean
\cite{Lau, Sangani1}. We call $\mathfrak{A}$
$\varphi$-\textit{amenable} if it is both left and right
$\varphi$-amenable.

Later in \cite{Sangani}, the authors introduced the concept of
$\varphi$-contractibility. Let $\mathfrak{A}$ be a Banach algebra
and $E$ be a Banach $\mathfrak{A}$-bimodule.  A continuous linear
operator $D:\mathfrak{A} \longrightarrow E$ is a \textit{derivation}
if it satisfies $ D(ab) = D(a)\cdot b + a \cdot D(b) $ for all $a,b
\in \mathfrak{A}$. Given $x \in E$, the \textit{inner} derivation
$ad_x : \mathfrak{A} \longrightarrow E$ is defined by $ad_x(a) =
a\cdot x - x\cdot a$. Let $ \varphi \in \Delta (\mathfrak{A}) $. We
write $ \mathbb{M}_\varphi$ $ [ _\varphi \mathbb{M}  ]$ for the set
of all Banach $\mathfrak{A}$-bimodules $E$ such that the right
[left] module action of $\mathfrak{A}$ on $E$ is given by $ x \cdot
a := \varphi(a) x$ $[ a \cdot x := \varphi(a) x]$ for $a \in
\mathfrak{A}$, $x \in E$. Precisely, $\mathfrak{A}$ is
\textit{right} [\textit{left}] $\varphi$-\textit{contractible} if
for each Banach $\mathfrak{A}$-bimodule $E \in \mathbb{M}_\varphi$ $
[E \in$ $_\varphi \mathbb{M}]$, every derivation $D:\mathfrak{A}
\longrightarrow E$ is inner. We say $\mathfrak{A}$ is
$\varphi$-contractible if it is both left and right
$\varphi$-contractible.

Recently and motivated by the above notions, several authors have
defined and studied the concept of $\varphi$-Connes amenability,
where $\varphi$ is a $w^*$-continuous homomorphism on a dual Banach
algebra \cite{Gh-Ja, mah, Sh-Am}. However, \cite{mah} was released
into the public domain over three years ago.

In this brief note, we take a look at $\varphi$-contractibility and
$\varphi$-Connes amenability in minute detail. While the authors put
a lot of effort into studying these concepts, their effort may have
been wasted here. Considering all three types of the notion of
$\varphi$-Connes amenability introduced in \cite{Gh-Ja, mah, Sh-Am},
we shall see that \textit{none} of them are new. In fact they
coincide with both $\varphi$-amenability and
$\varphi$-contractibility. Next, we shall prove that the concept of
$\varphi$-contractibility is also equivalent to an existing notion.
On closer inspection, saying that a Banach algebra $ \mathfrak{A}$
is $\varphi$-contractible is equivalent to saying that the
one-dimensional Banach $ \mathfrak{A}$-bimodule $\mathbb{C}_\varphi$
is projective. Although this concept goes back to Helemskii's works
in the 1970s (see his book \cite{Hel}, or alternatively a paper of
White \cite{Wh}), most of the authors who have studied
$\varphi$-contractibility seem unaware of this fact.

\section{$\varphi$-Connes amenability }

Let $\mathfrak{A}$ be a Banach algebra. A Banach
$\mathfrak{A}$-bimodule $E$ is \textit{dual} if there is a closed
submodule $E_*$ of $E^*$ such that $E = (E_*)^*$. We say $E_*$ the
\textit{predual} of $E$. A Banach algebra is \textit{dual} if it is
dual as a Banach $\mathfrak{A}$-bimodule. We write $\mathfrak{A} =
(\mathfrak{A}_*)^*$ if we wish to stress that $\mathfrak{A}$ is a
dual Banach algebra with predual $\mathfrak{A}_*$.

We start with the definition of $\varphi$-Connes amenability in the
sense of \cite{Gh-Ja}.
 Let $
\mathfrak{A}=( \mathfrak{A}_*)^*$ be a dual Banach algebra, and let
$\varphi \in \Delta (\mathfrak{A}) \cap \mathfrak{A}_*$. A dual
Banach $\mathfrak{A}$-bimodule $E \in$ $_\varphi \mathbb{M}  $ is
\textit{normal} if the module action $ a  \longmapsto x \cdot a$ of
$\mathfrak{A}$ on $E$ is $w^*$-continuous. A dual Banach algebra $
\mathfrak{A}=( \mathfrak{A}_*)^*$ is \textit{left}
$\varphi$-\textit{Connes amenable} if for every normal dual Banach $
\mathfrak{A}$-bimodule $E \in$ $_\varphi \mathbb{M}  $, every
$w^*$-continuous derivation $ D : \mathfrak{A} \longrightarrow E$ is
inner. Although they consider just left $\varphi$-Connes amenable
Banach algebras, there are similar definitions for right
$\varphi$-Connes amenable and $\varphi$-Connes amenable Banach
algebras. The authors show that (left) $\varphi$-Connes amenability
of $ \mathfrak{A}$ is equivalent to the existence a (left)
$\varphi$-mean \cite[Theorem 2.3]{Gh-Ja}.

\begin{theo} \label{2.1} Suppose that $\mathfrak{A}=( \mathfrak{A}_*)^*$ is a dual Banach
algebra and $\varphi \in \Delta (\mathfrak{A}) \cap \mathfrak{A}_*$.
Then the following statement are equivalent:

$(i)$ $\mathfrak{A}$ is $\varphi$-Connes amenable in the sense of
\cite{Gh-Ja};

$(ii)$ $\mathfrak{A}$ is  $\varphi$-contractible;

$(iii)$ $\mathfrak{A}$ is $\varphi$-amenable.
\end{theo}
{\bf Proof.}  The implications $(ii) \Longrightarrow (iii)$ and
$(iii) \Longrightarrow (i)$ are immediate.

$(i) \Longrightarrow (ii)$ Take a $\varphi$-mean $ m \in
\mathfrak{A}^{**}$. Consider the $\mathfrak{A}$-bimodule inclusion
map $ \imath : \mathfrak{A}_* \longrightarrow \mathfrak{A}^*$.
Taking adjoints, we obtain a $w^*$-$w^*$-continuous
$\mathfrak{A}$-bimodule map $ \xi :\mathfrak{A}^{**} \longrightarrow
\mathfrak{A}$. Now put $u = \xi(m) \in \mathfrak{A}$. It is easily
checked that $\varphi(u)=1$ and $ua = a u = \varphi(a) u$, for all
$a \in \mathfrak{A}$. Therefore by Theorem 3.1 below, $\mathfrak{A}$
is  $\varphi$-contractible.

$(ii) \Longrightarrow (i)$ Again by Theorem \ref{3.1}, there is an
element $ u \in \mathfrak{A}$ satisfying $\varphi(u)=1$ and $ua = ua
= \varphi(a) u$, for all $a \in \mathfrak{A}$. It is readily seen
that $u$ is a $\varphi$-mean and whence $\mathfrak{A}$ is
$\varphi$-Connes amenable. \qed

Now, we consider the definition of $\varphi$-Connes amenability in
the sense of \cite{Sh-Am}. Let $ \mathfrak{A}=( \mathfrak{A}_*)^*$
be a dual Banach algebra, and let $\varphi $ be a non-zero
$w^*$-continuous multiplicative linear functional on $
\mathfrak{A}$. The authors in \cite[Definition 2.1]{Sh-Am} say that
$ \mathfrak{A}$ is (\textit{left}) $\varphi$-\textit{Connes
amenable} if there exists $m \in \mathfrak{A}$ such that $
m(\varphi) = 1$ and $a m = \varphi(a) m$, for every $a \in
\mathfrak{A}$. Then by Theorem \ref {3.1}, left $\varphi$-Connes
amenability is nothing but right $\varphi$-contractibility. Hence
the following is straightforward.

\begin{theo} \label{4.1} Suppose that $\mathfrak{A}$ is a dual Banach
algebra and $\varphi $ be a non-zero $w^*$-continuous multiplicative
linear functional on $ \mathfrak{A}$. Then the following statement
are equivalent:

$(i)$ $\mathfrak{A}$ is $\varphi$-Connes amenable (in the sense of
\cite{Sh-Am});

$(ii)$ $\mathfrak{A}$ is  $\varphi$-contractible;

$(iii)$ $\mathfrak{A}$ is $\varphi$-amenable.
\end{theo}

Let $\mathfrak{A} $ be a dual Banach algebra and let $E$ be a Banach
$\mathfrak{A}$-bimodule. From \cite{Rund} we write $\sigma wc(E)$
for the set of all elements $ x \in E$ such that the maps
$$ \mathfrak{A} \longrightarrow E \ \ , \ \ \ a \longmapsto \left \{
\begin{array}{ll}
                        a \ . \ x \\
                        x \ . \ a
                      \end{array} \right. \ , $$ are $w^*$-weak
continuous.

 We conclude by looking at the definition of $\varphi$-Connes amenability from \cite{mah}. Suppose that
$\mathfrak{A}$ is a dual Banach algebra and $\varphi$ is a
homomorphism from $\mathfrak{A}$ onto $\mathbb{C}$. Then it is an
easy observation that $\varphi$ is $w^*$-continuous if and only if $
\varphi \in \sigma wc (\mathfrak{A}^*)$. Suppose that $\mathfrak{A}$
is a dual Banach algebra and $\varphi $ is a $w^*$-continuous
homomorphism from $\mathfrak{A}$ onto $\mathbb{C}$. We call
$\mathfrak{A}$ (\textit{right}) $\varphi$-\textit{Connes amenable}
if $\mathfrak{A}$ admits a (\textit{right}) $\varphi$-\textit{Connes
mean} $m$, i.e., there exists a bounded linear functional $m$ on $
\sigma wc (\mathfrak{A}^*)$ satisfying $ m(\varphi) = 1$ and $ m ( f
\ . \ a) = \varphi(a) m(f)$ for all $a \in \mathfrak{A}$ and $ f \in
\sigma wc (\mathfrak{A}^*)$. Similarly, we may consider
\textit{left} $\varphi$-Connes amenability. Meanwhile,
$\mathfrak{A}$ is $\varphi$-\textit{Connes amenable} if it is both
left and right $\varphi$-Connes amenable.

\begin{theo} \label{4.2} Suppose that $\mathfrak{A}=( \mathfrak{A}_*)^*$ is a dual Banach
algebra and $\varphi $ be a non-zero $w^*$-continuous multiplicative
linear functional on $ \mathfrak{A}$. Then the following statement
are equivalent:

$(i)$ $\mathfrak{A}$ is $\varphi$-Connes amenable (in the sense of
\cite{mah});

$(ii)$ $\mathfrak{A}$ is  $\varphi$-contractible;

$(iii)$ $\mathfrak{A}$ is $\varphi$-amenable.
\end{theo}
{\bf Proof.} Only $(i)\Longleftrightarrow (ii)$ needs the proof. Let
$ \iota : \mathfrak{A} \longrightarrow \sigma wc (\mathfrak{A}^*)^*$
be the $\mathfrak{A}$-bimodule map obtained by composing the
canonical inclusion $\mathfrak{A} \longrightarrow \mathfrak{A}^{**}$
with the quotient map $ \mathfrak{A}^{**} \longrightarrow \sigma wc
(\mathfrak{A}^*)^*$, so that $ \langle \iota(a) , \psi \rangle =
\psi (a)$ for all $a \in \mathfrak{A}$ and $ \psi \in \sigma wc
(\mathfrak{A}^*)$.

$(i) \Longrightarrow (ii)$ Since $\mathfrak{A}$ is a dual Banach
algebra, $\mathfrak{A}_*$ is an $\mathfrak{A}$-bimodule and
$\mathfrak{A}_* \subseteq \sigma wc (\mathfrak{A}^*)$
\cite[Corollary 4.6]{Rund}. Therefore taking adjoints gives us a
$w^*$-$w^*$-continuous $\mathfrak{A}$-bimodule map $ \xi : \sigma wc
(\mathfrak{A}^*)^* \longrightarrow \mathfrak{A}$. Notice that $\xi
\circ \iota (a) = a$ for all $a \in \cA$. By the assumption, there
exists a $\varphi$-Connes mean $m \in \sigma wc (\mathfrak{A}^*)^*$.
Setting $u = \xi(m) \in \mathfrak{A}$, we observe that
$\mathfrak{A}$ is $\varphi$-contractible by Theorem \ref{3.1}.

$(ii) \Longrightarrow (i)$ Take $ u \in \mathfrak{A}$ satisfying
$\varphi(u)=1$ and $u a = au = \varphi(a) u$, for all $a \in
\mathfrak{A}$. Then $\iota(u)$ is a $\varphi$-Connes mean on $
\sigma wc (\mathfrak{A}^*)$. \qed

\section{$\varphi$-contractibility }

 It was shown that right [left] $\varphi$-contractibility of
$\mathfrak{A}$ is equivalent to the existence of a \textit{right}
[\textit{left}] $\varphi$-\textit{diagonal} for $\mathfrak{A}$, i.
e., an element $ m \in \mathfrak{A} \widehat{\otimes} \mathfrak{A}$
such that $ \varphi ( \pi (m) ) = 1$ and $ a\cdot m = \varphi (a) m
$ $[m \cdot a =  \varphi (a) m]$ for $a \in \mathfrak{A}$, where $
\pi : \mathfrak{A} \widehat{\otimes} \mathfrak{A} \longrightarrow
\mathfrak{A}$ is the bounded linear map determined by $ \pi (a
\otimes b ) = ab$. If $m$ is both left and right $\varphi$-diagonal,
it called $\varphi$-\textit{diagonal}.

The following is likely to be well-known, but since we could not
locate a reference, we include a proof.

\begin{theo} \label{3.1}Suppose that $\mathfrak{A}$ is a Banach
algebra and $\varphi \in \Delta(\mathfrak{A} ) $. Then
$\mathfrak{A}$ is $\varphi$-contractible if and only if there exists
an element $u \in \mathfrak{A}$ satisfying $$  (*) \ \ \ \ \ \ \ \ \
\ \ \ \ \ \ \varphi(u) = 1 \ \ \text{and} \ \ u a = a u = \varphi(a)
u \ \ \ \ \ \ \ (a \in \mathfrak{A}) \ .$$
\end{theo}
{\bf Proof.} Let $u \in \mathfrak{A}$ satisfies the conditions in
$(*)$. Let $ D : \mathfrak{A} \longrightarrow E$ be a derivation for
a Banach $\mathfrak{A}$-bimodule $E \in \mathbb{M}_\varphi$. It is
routinely checked that $ D^* ( f \ . \ a ) = D^* ( f ) \ . \ a  -
\langle
  Da , f \rangle \varphi$, for each $f \in E^*$ and $a \in
  \mathfrak{A}$. Put $ t := Du \in E$. For $f \in E^*$ and $a \in
  \mathfrak{A}$ we have
\begin{align*}
\langle  f, a \ . \ t \rangle &= \langle D^* ( f \ . \ a ) , u
\rangle = \langle  D^* ( f ) , a u \rangle - \langle Da , f \rangle
\\&= \varphi(a) \langle  D^* ( f ) , u \rangle - \langle Da , f
\rangle = \varphi(a) \langle  f , t \rangle - \langle Da , f \rangle
\ .
\end{align*}
Therefore $a \ . \ t = \varphi(a) t - Da $, $a \in
  \mathfrak{A}$, and hence $ D = ad_{-t}$. Whence $\mathfrak{A}$ is
  right $\varphi$-contractible. A similar argument shows that $\mathfrak{A}$
  is also left $\varphi$-contractible.

Conversely, suppose that $ m \in \mathfrak{A} \widehat{\otimes}
\mathfrak{A}$ is a $\varphi$-diagonal for $\mathfrak{A}$. Put $u:=
\pi(m) \in \mathfrak{A}$. Now, it is easily checked that $u$ has the
desired properties in $(*)$.\qed

Let $\mathfrak{A}$ be a Banach algebra with the \textit{unitization}
$\mathfrak{A}^\sharp $, and let $P$ be a left Banach
$\mathfrak{A}$-module. We recall that $P$ is a \textit{projective}
left $\mathfrak{A}$-module if the multiplication map $$   \pi :
\mathfrak{A}^\sharp \widehat{\otimes} P \longrightarrow P \ ; \ \ a
\otimes x \longmapsto a \cdot x  \ \ \ (a \in  \mathfrak{A}^\sharp,
x \in P)$$ has a right inverse which is also a left
$\mathfrak{A}$-module homomorphism. Similar definitions hold for
projective right $\mathfrak{A}$-modules and projective
$\mathfrak{A}$-bimodules.

For $\varphi \in \Delta(\mathfrak{A} ) $, the space $
\mathbb{C}_\varphi = \{ \alpha \varphi \ : \ \alpha \in \mathbb{C}
\}$ is a Banach $\mathfrak{A}$-bimodule with module actions $ a
\cdot \varphi = \varphi \cdot a := \varphi(a) \varphi$, $(a \in
\mathfrak{A})$.

\begin{theo} \label{3.2}Suppose that $\mathfrak{A}$ is a Banach
algebra and $\varphi \in \Delta(\mathfrak{A} ) $. Then
$\mathfrak{A}$ is $\varphi$-contractible if and only if $
\mathbb{C}_\varphi$ is a projective Banach $\mathfrak{A}$-bimodule.
\end{theo}
{\bf Proof.} Without loss of generality, we may assume that
$\mathfrak{A}$ is unital. Let $ \mathbb{C}_\varphi$ be projective as
a left $\mathfrak{A}$-module. Then there exists a bounded linear map
$ \rho : \mathbb{C}_\varphi \longrightarrow \mathfrak{A}
\widehat{\otimes} \mathbb{C}_\varphi$ satisfying $ \pi \rho
(\varphi) = \varphi$ and $a \cdot \rho (\varphi) = \varphi(a) \rho
(\varphi)$ for each $a \in \mathfrak{A}$. We have $ \rho (\varphi) =
\sum_{n=1}^\infty a_n \otimes \varphi$, where $a_n \in \mathfrak{A}$
$(n=1,2,...)$ with $ \sum_{n=1}^\infty || a_n || < \infty $. Putting
$ u:=\sum_{n=1}^\infty a_n \in \mathfrak{A}$, we observe that
$\varphi (u) = 1$ and $ a u = \varphi(a) u$, $a \in \mathfrak{A}$.
Now, by Theorem \ref{3.1} $\mathfrak{A}$ is right
$\varphi$-contractible.

Conversely, let $\mathfrak{A}$ be right $\varphi$-contractible. Take
$ u \in \mathfrak{A}$ with $\varphi (u) = 1$ and $ a u = \varphi(a)
u$ for all $a \in \mathfrak{A}$. Then it is easy to verify that the
map $ \rho : \mathbb{C}_\varphi \longrightarrow \mathfrak{A}
\widehat{\otimes} \mathbb{C}_\varphi$ defined by $ \rho (\varphi) :=
u \otimes \varphi$ is a left $\mathfrak{A}$-module homomorphism
which is a right inverse of $\pi$. Whence $\mathbb{C}_\varphi$ a is
projective left $\mathfrak{A}$-module.

Similarly, one can see that $\mathbb{C}_\varphi$ is a projective
right $\mathfrak{A}$-module if and only if $\mathfrak{A}$ is left
$\varphi$-contractible. \qed

\end{document}